\input amstex
\documentstyle{amsppt}

\topmatter

\title
Rational Obstruction Theory and Rational Homotopy Sets
\endtitle

\author Martin Arkowitz and Gregory Lupton \endauthor

\abstract
We develop an obstruction theory for homotopy of
homomorphisms $f,g : {\Cal M }\to{\Cal N }$ between minimal
differential graded algebras.
We assume that ${\Cal M }=\Lambda V$ has an obstruction
decomposition given by
$V=V_0\oplus V_1$ and that $f$ and $g$ are homotopic on $\Lambda
V_0$.  An obstruction  is then obtained as a vector space
homomorphism $V_1\to H^*({\Cal N})$.   We investigate the
relationship between the condition that $f$ and $g$ are homotopic
and the  condition that the obstruction is zero.
The obstruction theory is then applied to study
the set of homotopy classes  $[{\Cal M },{\Cal N }]$.
This enables us to give a fairly complete answer to a conjecture of
Copeland-Shar on the size of the homotopy set $[A,B]$ when $A$ and
$B$ are rational spaces.  In addition, we give examples of minimal
algebras (and hence of rational spaces) that have few homotopy
classes of self-maps.
\endabstract



\address
Department of Mathematics, Dartmouth College, Hanover NH 03755
\endaddress
\email
Martin.Arkowitz\@dartmouth.edu
\endemail

\address
Department of Mathematics, Cleveland State University, Cleveland OH
44115
\endaddress
\email
Lupton\@math.csuohio.edu
\endemail

\endtopmatter

\document

\subhead{\bf Introduction}\endsubhead   A basic object of study in
unstable homotopy theory, perhaps {\it the} basic object of study,
is the set of homotopy classes of maps from one topological space
to another. For arbitrary spaces, it is difficult to describe this
homotopy set fully, and the best that can be  hoped for is some
partial information.  In this paper we use methods from rational
homotopy theory to study this  set for rational spaces (i.e.,
spaces whose homotopy groups are rational vector spaces). This can
be regarded as an approximation to homotopy sets for finite
complexes.

As is well known, the homotopy theory of rational spaces is
equivalent to the homotopy theory of minimal, differential, graded
commutative algebras over the rationals (minimal algebras, for
short).   Minimal algebras provide an effective algebraic setting
to work in, and we begin by considering homomorphisms of minimal
algebras. We develop an obstruction theory for homotopy of
homomorphisms of minimal algebras and present a few applications
of this  theory.  These include  a simple proof that there are no
non-trivial phantom maps between minimal algebras and results on a
conjecture of Copeland-Shar that the  homotopy set $[{\Cal M},
{\Cal N}]$ for minimal algebras ${\Cal M}$ and ${\Cal N}$ is
either trivial or infinite. In addition, we give examples of
minimal algebras that have few self-maps,  including an elliptic
minimal algebra with trivial group of homotopy classes of
self-equivalences.

Because of the categorical equivalence mentioned above, the
results on homotopy classes of homomorphisms of minimal algebras,
obtained by the  obstruction theory, translate immediately into
corresponding results about homotopy classes of maps of rational
spaces. Thus we obtain results, of interest in their own right,
about the set of homotopy classes of maps $[A, B]$, where $A$ and
$B$ are rational spaces.    For instance, in this way we obtain an
example of a non-trivial rational space with trivial group of
self-equivalences.

Next we describe the  organization of the paper.  Section 1
establishes our notation and basic conventions.  The obstruction
theory is developed in Section 2.  We consider maps of minimal
algebras $f \: {\Cal M} \to {\Cal N}$, where ${\Cal M} = \Lambda
V$ and $V$ admits a certain direct sum decomposition $V = V_0
\oplus V_1$. Then if  $f, g \: \Lambda V \to {\Cal N}$ are two
maps whose restrictions to $\Lambda V_0$ are homotopic, we define
an obstruction homomorphism $V_1 \to H^*({\Cal N})$. Theorem 2.5
asserts that $f$ and $g$ are homotopic via a homotopy that extends
the given homotopy if and only  if the obstruction is zero.
Therefore if the obstruction is zero, the maps  are homotopic.
However, the converse is not true in general. In Propositions 2.6
and 2.8 we establish the converse in special cases  which suffice
for our applications.

In Section 3 the obstruction theory is generalized and we obtain
the result about phantom maps, mentioned above (Theorem 3.1).
Also, in Corollary 3.3, we give a  necessary condition for a map
of arbitrary  minimal algebras to be homotopically non-trivial.
This is later used  to show a certain homotopy set is infinite.

In Section 4 we investigate the following conjecture:  If the set
of homotopy classes of maps between two rational spaces is finite,
then it consists of one element.  In Theorem 4.6 we prove that
this is true  if either of the spaces belongs to a very broad
class of spaces, namely, universal spaces. This is done by
translating the conjecture into one about minimal algebras and
then applying the previous obstruction theory. However, in Section
5 we also give examples of non-universal, rational spaces that
have exactly two and exactly three homotopy classes of self-maps.
These examples have several interesting features such as having
groups of self-homotopy equivalences the trivial group and ${\Bbb
Z}_2$, respectively.

\subhead{\bf \S1 Preliminaries}\endsubhead In general,   our
notation and conventions follow the references \cite{G-M},
\cite{Ha$_2$}, \cite{H-S} and \cite{Su$_2$}.   By a vector space
we mean a graded vector space over the  rational numbers ${\Bbb
Q}$, i.e., a collection $V=\{ V^k \ |\ k\ \hbox {an integer} \ge 0
\}$, such that each $V^k$ is a  vector space over ${\Bbb Q}$.   If
$v_1,\ldots ,v_r,\ldots $ is a basis of $V$, that is, $v_1,\ldots
,v_{i_0}$ is a basis of $V^0$, $v_{i_0+1},\ldots ,v_{i_1}$ is a
basis of $V^1$, etc., then we write $V=\langle v_1, \ldots
,v_r,\ldots \rangle $. If the set of basis vectors of $V$ is
finite, we say that $V$ is finite-dimensional.

We consider differential graded commutative algebras $({\Cal
A},d)$ over ${\Bbb Q}$ --- called DG algebras --- with
differential $d$ a derivation of degree $+1$.  We write $x\in
{\Cal A}$ to indicate that $x \in {\Cal A}^n$ for some $n\geq0$
and let $|x| = n$ be the degree of $x$.   We denote the cohomology
algebra of ${\Cal A}$ by $H^*({\Cal A})$ and let $[x]\in H^*({\Cal
A})$ stand for the cohomology class of the cocycle $x \in {\Cal
A}$.     By a map $f\:{\Cal A} \to {\Cal B}$ of DG algebras we
mean a DG algebra homomorphism.  Then $f$ induces a map of
cohomology algebras $f^*\: H^*({\Cal A}) \to H^*({\Cal B})$.  The
identity map of ${\Cal A}$  will  be denoted by $\iota \: {\Cal A}
\to {\Cal A}$ and the zero map by $0 \: {\Cal A} \to {\Cal B}$.

The free graded commutative algebra generated by $V$ is denoted
$\Lambda V$ and a basis for $V$ is  called a set of algebra
generators for $\Lambda V$. If $V=\langle v_1, \ldots ,v_r,\ldots
\rangle $, we write $\Lambda V = \Lambda ( v_1, \ldots
,v_r,\ldots )$.  A DG algebra $({\Cal M},d)$ is a {\it minimal
algebra} if \ (i) ${\Cal M} =\Lambda V$ for some vector space $V$
and \ (ii) the image of $d$ is decomposable, i.e., a sum of
products of positive dimensional elements (cf. \cite{G-M}). A {\it
cofibration\/}
 is a map
$\Lambda V \to \Lambda (V \oplus U)$ of minimal algebras such that
$U$ has a basis $u_1,\ldots ,u_r,\ldots$ that satisfies $d(u_r)
\in \Lambda(V \oplus\langle u_1,\ldots ,u_{r-1}\rangle)$ for each
$r$ and $|u_1|\le |u_2| \le \cdots $ \cite{Ha$_2$, p.\,14}.

We will assume that the following connectedness and finiteness
conditions always hold.
 If ${\Cal A}$ is a DG algebra, then
$H^0({\Cal A}) = {\Bbb Q}$, $H^1({\Cal A})=0$ and $H^k({\Cal A})$
is finite-dimensional for each $k$. If ${\Cal M}$ is a minimal
algebra, then ${\Cal M}^0 = {\Bbb Q}$, ${\Cal M}^1=0$ and ${\Cal
M}^k$ is finite-dimensional for every $k$.  Thus if ${\Cal M}
=\Lambda V$, then $V^0 =V^1 =0$ and each $V^k$ is
finite-dimensional.

If ${\Cal A} = \Lambda V$ is a free algebra, then we specify a
{\it second grading} $\{ {\Cal A}_r |\, r\ge 0\}$ on ${\Cal A}$ by
giving a basis of $V$ a second grading and  extending
multiplicatively. We then denote $\oplus_{r \leq n} V_r$ by
$V_{(n)}$, and  the free algebra $\Lambda(V_{(n)})$ by $\Lambda
V_{(n)}$.  If ${\Cal A}$ has a second grading, then we  refer to
${\Cal A}$ as a {\it bigraded algebra}.  A {\it bigraded DG
algebra\/} is a bigraded algebra ${\Cal A}$  such that $d\: {\Cal
A}_{r} \to {\Cal A}_{r-1}$.     If ${\Cal A}$ is any bigraded DG
algebra, then the second grading on ${\Cal A}$ makes the
cohomology algebra $H^*({\Cal A})$ into a bigraded algebra. An
example of a bigraded minimal algebra is the bigraded model of a
formal space, introduced in \cite{H-S}. More generally, a minimal
algebra ${\Cal M} =\Lambda V$ is a {\it filtered} minimal algebra
if ${\Cal M}$ is bigraded as an algebra, $V=\oplus_{i \geq 0} V_i$
and $d|_{V_n} \: V_n \to \Lambda V_{(n-1)}$.    In this case there
is a sequence of inclusions $\Lambda V_{0} \subseteq \cdots
\subseteq \Lambda V_{(r)} \subseteq \Lambda V_{(r+1)} \subseteq
\cdots$  and ${\Cal M} = \bigcup_r \Lambda V_{(r)}$.

For maps of DG algebras we use the notion of homotopy described in
\cite{Ha$_2$, Ch.\,5}, \cite{H-S, p.\,240} and \cite{Su$_2$, \S3}
which we now review.  Suppose $\Lambda V$ is a minimal algebra.
Define the DG algebra $\Lambda V^I =\Lambda(V\oplus \overline {V}
\oplus \widehat V)$ with  differential  $d$ as follows: $\widehat
V$ is an isomorphic copy of $V$ and $\overline {V}$ is the
desuspension of $V$, i.e., $\overline {V} ^p=V^{p+1}$.
Furthermore,
 $d$  agrees with the differential on
$\Lambda V$, $d(\bar v) = \hat v$ and $d(\hat v)=0$, for $\bar v
\in \overline {V}$ and $\hat v \in \widehat V$ corresponding to $v
\in V$.  In addition, there is a degree $-1$ derivation $i\:
\Lambda V^I \to \Lambda V^I$ defined on generators by $i(v) =\bar
v$, $i(\bar v ) = 0$ and $i(\hat v) = 0$.  We then obtain a degree
0 derivation $\gamma \: \Lambda V^I \to \Lambda V^I$ by setting
$\gamma =di+id$. Finally, we have a map $\alpha \: \Lambda V^I \to
\Lambda V^I$ given by
$$\alpha = \sum^{\infty}_{n=0} {1\over {n!}}\gamma ^n.$$
If $f,g \: \Lambda V \to {\Cal N}$ are maps, then  a {\it homotopy
H from f to g\/} (or   {\it a homotopy H which starts at f and
ends at g \/}) is a map $H \: \Lambda V^I \to {\Cal N}$ such that
$H|_{\Lambda V} =f$ and $H\alpha |_{\Lambda V} =g$.   We then
write $f \simeq_H g$ or just $f \simeq g$.  This is an equivalence
relation on the set of maps between two minimal algebras.    We
denote the collection of homotopy classes of maps from ${\Cal M}$
to ${\Cal N}$ by $[{\Cal M},{\Cal N}]$ with the  homotopy class of
$f\: {\Cal M} \to {\Cal N}$ written $[f]$.  A map which is
homotopic to the zero map is said to be nullhomotopic.

Let $\Lambda V$  and $\Lambda W$ be minimal algebras and $j \:
\Lambda V \to \Lambda W$  a map.   Then there is a DG map $j^I \:
\Lambda V^I \to \Lambda W^I$  extending $j$ defined as follows:
 $j^I(\bar v) = i j(v)$ and
$j^I(\hat v) = di j(v)$ for $\bar v \in \overline V$ and $\hat v
\in \widehat V$. Then  we say that $j \: \Lambda V \to \Lambda W$
has the {\it homotopy extension property\/} if, whenever $f \:
\Lambda W \to {\Cal N}$ is a map and $H \: \Lambda V^I\to {\Cal
N}$ is a homotopy which  starts at $fj$,  there is a homotopy
${\widetilde H} \: \Lambda W^I \to {\Cal N}$  which starts at $f$
such that ${\widetilde H}j^I = H$. The following lemma is clear.

\proclaim{\bf 1.1 Lemma}  A cofibration $j \: \Lambda V \to
\Lambda (V \oplus U)$ has the homotopy extension property.
\endproclaim

Finally, we say that a topological space is  of {\it finite
type\/} if its rational homology is finite-dimensional in each
degree.  All topological spaces in this paper are based and have
the based homotopy type of a 1-connected CW-complex of finite
type. For a space $X$, we denote by $X_{\Bbb Q}$ the
rationalization of $X$, and for a map $\phi\: X \to Y$ of spaces,
we denote by $\phi_{\Bbb Q}$ the rationalization of $\phi$
\cite{H-M-R}.   Then $\phi $ is called a ${\Bbb Q}$-equivalence if
$\phi _{\Bbb Q}$ is a homotopy equivalence.

\subhead{\bf \S2 Obstruction Theory}\endsubhead In this section we
briefly sketch an obstruction theory for homotopy of maps of
minimal algebras $f,g \: {\Cal M} \to {\Cal N}$.

\proclaim {2.1 Definition}  Let ${\Cal M} = \Lambda V$ be a
minimal algebra.  An {\it obstruction decomposition\/} for ${\Cal
M}$ consists of a decomposition $V = V_0\oplus V_1$ of $V$ as a
direct sum of (graded) vector spaces $V_0$ and $V_1$, such that
$d|_{V} \: V \to \Lambda V_0$.
\endproclaim

Note that if $V = V_0\oplus V_1$ is an obstruction decomposition
for $\Lambda V$, then $\Lambda V_0 \to \Lambda(V_0\oplus V_1)$ is
a cofibration. We do not require $d(V_0) = 0$. We now describe the
 main examples of an obstruction decomposition.

\smallskip

\noindent{\bf 2.2 Examples.} (1) If  $V$ is finite-dimensional, $V
=  V^2 \oplus \cdots \oplus V^n$ for some $n$.   Set $V_0 =
V^{(n-1)} = V^2 \oplus \cdots \oplus V^{n-1}$ and $V_1 = V^n$.  We
call this obstruction decomposition the {\it degree
decomposition}.

\noindent{}(2) If ${\Cal M}$ has an obstruction decomposition $V =
V_0 \oplus V_1$ with $d(V_0)=0$,   then we say that the
obstruction decomposition is a {\it two-stage decomposition}, and
that ${\Cal M}$ is a {\it two-stage minimal algebra}.  Two-stage
minimal algebras are considered in \cite{H-T}, where it is  shown
that the Sullivan minimal model of every homogeneous space  is a
two-stage minimal algebra \cite{H-T, \S4}.

\noindent{}(3) Let ${\Cal M} = \Lambda V$ be a filtered minimal
algebra as in Section 1 with $V = \oplus_{i \geq 0} V_i$.  Then
for each fixed $n \geq 1$, $\Lambda V_{(n)}$ has an obstruction
decomposition $V = U_0\oplus U_1$ where $U_0 = V_{(n-1)}$ and $U_1
= V_n$.  This arises, for example, when ${\Cal M}$ is formal and
the filtration of ${\Cal M}$ is given by the bigraded model.

\smallskip

We assume throughout this section that ${\Cal M} =\Lambda V$ is a
minimal algebra with an obstruction decomposition $V = V_0 \oplus
V_1$.    Since $\Lambda V_0$ is a minimal algebra contained in
$\Lambda V$, the DG algebra $\Lambda V_0^I$ is a sub-DG algebra of
$\Lambda V^I$. If $S \subseteq \Lambda V_0^I$ is a set, then $( S
) \subseteq \Lambda V_0^I$  denotes the ideal of $\Lambda V_0^I$
generated by $S$.

The following lemma plays a key role in this section.

\proclaim{\bf 2.3 Lemma}  If $x \in V$, then $\alpha (x) = x +\hat
x + \xi$, where $\xi$ is a decomposable element in the
intersection of ideals $\big( \overline{V_0} \big) \cap \big( V_0
\oplus \widehat{V_0} \big)$ of $\Lambda V_0^I$. If $d(x) = 0$,
then $\alpha (x) = x + \hat x$.
\endproclaim

\demo{Proof}  Since $\gamma(x) = di(x) + id(x) = \hat{x} + id(x)$
and $\gamma(\hat{x}) = 0$, we have $\alpha(x) = x + \hat{x} +
\xi$, where
$$\xi = id(x) + {1 \over 2!}\gamma(id(x)) + \cdots + {1 \over
n!}\gamma^{n-1}(id(x)) + \cdots\ .$$
If $d(x) = 0$, then clearly $\alpha (x) = x + \hat x$.

Now $d(x)$ is decomposable and therefore so  is $\xi$.
 From the Definition 2.1,
$d(x) \in \Lambda V_0$.  Since  $\Lambda V_0^I$ is stable under
$i$ and $d$, and hence under $\gamma$, it follows that $\xi \in
\Lambda V_0^I$.

Next $\gamma( \overline{V} ) = 0$, so the ideal $\big(
\overline{V_0} \big)$  is stable under $\gamma$.  Since $id(x) \in
\overline{V_0}\cdot \Lambda^+ V_0$, where $\Lambda^+ V_0$ denotes
the positive-dimensional elements of $\Lambda V_0$, we have $\xi
\in \big( \overline{V_0} \big)$.

If $u \in V_0$, then $\gamma(u) = \hat u + id(u)$ with $id(u) \in
\overline{V_0}\cdot \Lambda^+ V_0$.  Since $\gamma(\widehat V) =
0$, it follows that for  $y \in V_0 \oplus \widehat{ V_0 }$, we
have $\gamma(y) \in \big( V_0 \oplus \widehat{ V_0 } \big)$, and
hence
 $\big( V_0 \oplus \widehat{ V_0 } \big)$  is
stable under $\gamma$.  Finally,  $id(x) \in \overline{V_0}\cdot
\Lambda^+ V_0 \subseteq  \big( V_0 \oplus \widehat{ V_0 } \big)$.
Therefore $\xi \in \big( V_0 \oplus \widehat{ V_0 } \big)$. \hfill
$\square$
\enddemo

Now consider a second minimal algebra ${\Cal N}$.   Suppose that
$f,g \: \Lambda V \to {\Cal N}$ are maps such that $f|_{\Lambda
V_0}  \simeq  g|_{\Lambda V_0}$ by a homotopy $H \: \Lambda V_0^I
\to {\Cal N}$.

\proclaim{2.4 Definition} The obstruction to $f$ and $g$ being
homotopic is the map ${\Cal O}^H_{f,g} \: V_1 \to H^*({\Cal N})$
of vector spaces defined for $w \in V_1$ by
$${\Cal O}^H_{f,g}(w) = [f(w) + H\big(\alpha(w) - w - \hat{w}\big)
- g(w)].$$

\endproclaim

Note that $\alpha(w) - w - \hat{w} \in \Lambda V_0^I$, so
$H(\alpha(w) - w - \hat{w})$ makes sense, and that $f(w) +
H\big(\alpha(w) - w - \hat{w}\big) - g(w)$ is a cocycle.

We now give a basic result of this section.

\proclaim{\bf 2.5 Proposition} Let $f,g \: \Lambda V \to {\Cal N}$
be maps such that $f|_{\Lambda V_0}  \simeq_H  g|_{\Lambda V_0}$.
Then ${\Cal O}^H_{f,g} = 0$ if and only if $f \simeq_K g$ for some
homotopy $K \: \Lambda V^I \to {\Cal N}$ which is an extension of
$H$.
\endproclaim

The proof is straightforward, and hence omitted.

\smallskip

It follows from Proposition 2.5, that if the obstruction ${\Cal
O}^H_{f,g}$ is zero, then $f$ and $g$ are homotopic.  However,  if
$f$ and $g$ are homotopic  by some homotopy $K$ that is not an
extension of $H$, then the obstruction ${\Cal O}^H_{f,g}$ is not
necessarily zero.  Indeed, a simple example in \cite{A-L,
Ex.\,4.6} illustrates this.

\smallskip

For our purposes, it is useful to have conditions under which $f
\simeq_K g$ for any homotopy $K$ implies that ${\Cal O}^H_{f,g}$
is zero.    We present two such results.

\proclaim{\bf 2.6 Proposition}  Let $f,g \: {\Cal M} =\Lambda V
\to {\Cal N}$ be  such that $f|_{\Lambda V_0}  =  g|_{\Lambda V_0}
= 0$ and let ${\Cal O}_{f,g}$ be the obstruction using the zero
homotopy. Then ${\Cal O}_{f,g} = {\Cal O}_{f,g}^H \: V_1 \to
H^*({\Cal N})$ for any homotopy $H$ from $f|_{\Lambda V_0}$ to
$g|_{\Lambda V_0}$, and so $ {\Cal O}_{f,g}^H(w) = [f(w)-g(w)]$,
for $w\in V_1$. Furthermore, ${\Cal O}_{f,g} = 0
\Longleftrightarrow
 f \simeq g \: \Lambda V \to {\Cal N}$.
\endproclaim

\demo{Proof} First we show that ${\Cal O}_{f,g} = {\Cal
O}_{f,g}^H$. If $w \in V_1$,  then ${\Cal O}_{f,g}(w) = [f(w) -
g(w)]$. For any homotopy $H\: \Lambda V_0^I \to {\Cal N}$ from
$f|_{\Lambda V_0} = 0$ to $g|_{\Lambda V_0} = 0$, we have ${\Cal
O}_{f,g}^H(w) = [f(w) + H(\xi) - g(w)]$ for some $\xi \in \big(
V_0 \oplus \widehat{V_0} \big)$ by Lemma 2.3.  The result will
follow by showing that $H\big(\big( V_0 \oplus \widehat{V_0}
\big)\big) = 0$. Clearly $H(v) = f(v) = 0$ for $v \in V_0$. We
claim that $H(\hat{v}) = 0$ and argue by induction on the degree
of $\hat{v}$. If $|\hat{v}|=|v| = 2$, then $d(v) = 0$ and
$\alpha(v) = v + \hat{v}$ by Lemma 2.3.  Hence $0 = g(v) =
H\alpha(v) = H(v) + H(\hat{v}) = H(\hat{v})$.  Now suppose that
$H(\hat{x}) = 0$ for all $x$ with $|x| < n$, and suppose $v \in
V_0$ with $|\hat{v}|=|v| = n$. Then $0 = g(v) = H\alpha(v) = H(v +
\hat{v} + \xi)$, with $\xi$ decomposable and in $\big( V_0 \oplus
\widehat{V_0} \big) $ by Lemma 2.3.  Thus $0 = H(v) + H(\hat{v}) +
H(\xi) = H(\hat{v}) + H(\xi)$. This gives $H(\hat{v}) = - H(\xi)$.
But $\xi$ is a decomposable element, each monomial of which
contains as a factor either some $x \in V$ or some $\hat{x} \in
\widehat{V}$ of degree $< n$.  Thus $H(\xi) = 0$ and so
$H(\hat{v}) = 0$.  This completes the induction, and so for all $v
\in V_0$, both $H(v)$ and $H(\hat{v})$ are zero. Therefore
$H\big(\big( V_0 \oplus \widehat{V_0} \big) \big) = 0$.

Now we prove the last assertion. Let $f \simeq _K g :\Lambda V \to
{\Cal N}$ for some homotopy $K$ and let $K'=K|_{\Lambda V_0^I}$.
Then by Proposition 2.5 and the first assertion of 2.6, ${\Cal
O}_{f,g} = {\Cal O}_{f,g}^{K'} =0$. \hfill$\square$
\enddemo

\noindent{\bf 2.7 Remark.}  Under the hypothesis of Proposition
2.6 we will use the notation ${\Cal O}_{f,g}$ to denote ${\Cal
O}_{f,g}^H$ for any homotopy $H$.  Furthermore, we  have zero
homotopies from $f|_{\Lambda V_0}$ to $0$ and from $g|_{\Lambda
V_0}$ to $0$ and so can write ${\Cal O}_{f,g} = {\Cal O}_{f,0} -
{\Cal O}_{g,0}$.
\smallskip

\proclaim{\bf 2.8 Proposition}  Let $f \: {\Cal M} =\Lambda V \to
{\Cal N}$ be   such that $f|_{\Lambda V_0} \simeq_H 0$. Then
 ${\Cal O}_{f,0}^H = 0 \Longleftrightarrow
f \simeq 0 \: {\Cal M} \to {\Cal N}$.
 \endproclaim

\demo{Proof}
  Suppose
$fj\simeq _{H} 0: \Lambda V_0 \to {\Cal N}$.  Since $j \: \Lambda
V_0 \to \Lambda V$ has the homotopy extension property by Lemma
1.1, there is a homotopy $G \: \Lambda V^I \to {\Cal N}$ that
starts at $f$ and extends $H$.  If ${\Cal O}_{f,0}^H = 0$, then $f
\simeq 0 $ by Proposition 2.5.  Conversely, suppose $f \simeq 0 $
by some homotopy. Since $G$ ends at $f' \: {\Cal M} \to {\Cal N}$,
then $f'|_{\Lambda V_0}=0$.  But $f' \simeq f\simeq 0$.    Thus by
Proposition 2.6 we have $0={\Cal O}_{f',0}$. Hence
$
0={\Cal O}_{f,f'}^{H} = {\Cal O}_{f,0}^{H}-{\Cal O}_{f',0},$ and
so ${\Cal O}_{f,0}^{H}=0$. \hfill $\square$
\enddemo

\subhead{\bf \S3 Extensions of the Obstruction Theory}\endsubhead
In this section we study homotopy of  homomorphisms of minimal
algebras $f,g \: {\Cal M }\to{\Cal N }$, where ${\Cal M }$ does
not necessarily admit an obstruction decomposition. We consider
when a map $f \: {\Cal M }\to{\Cal N }$ is nullhomotopic. This
will be needed in the sequel.

Suppose that ${\Cal M} = \Lambda V$ is a filtered minimal algebra
with $V =\oplus_{i\geq 0} V_i$ and $d|_{V_n} \: V_n \to \Lambda
V_{(n-1)}$ for each $n$,  where $V_{(n)} =\oplus _{0\leq i \leq
n}V_i$ (see Section 1).   Every minimal algebra has at least one
such filtration since we can always set $V_j = V^j$.  If ${\Cal
M}$ is formal, then the bigraded model of ${\Cal M}$ provides
another filtration. We denote the inclusion $j_n \: \Lambda
V_{(n)} \to \Lambda V$. If $f \: {\Cal M} \to {\Cal N}$ is any
map, then $f j_n \: \Lambda V_{(n)} \to {\Cal M } \to {\Cal N }$
will be denoted $f_n$.   Note that $\Lambda V_{(n)}$ admits an
obstruction decomposition $V_{(n)} \cong V_{(n-1)} \oplus V_n$ and
so the results of Section 2 apply to $\Lambda V_{(n)}$.

Now suppose we have homotopic maps $f \simeq _H g \: {\Cal M }  =
\Lambda V \to {\Cal N }$. By composing the homotopy $H$ with
$j_n^I$, we have $f_n \simeq  g_n \: \Lambda V_{(n)} \to {\Cal N
}$, for each $n$. The next result  is a converse of this in a
special case.

\proclaim{3.1 Theorem}  Suppose $f \: \Lambda V \to {\Cal N}$ is a
map of minimal algebras and ${\Cal M}$ admits a filtration given
by $V = \oplus_{i\geq 0} V_i$.  If $f_n \: \Lambda V_{(n)} \to
{\Cal N }$ is nullhomotopic for each $n$,  then $f \: \Lambda V
\to {\Cal N }$ is nullhomotopic.
\endproclaim

\demo{Proof} We argue inductively over $n$ to show that there is a
sequence of homotopies $H_n \: \Lambda V_{(n)}^I \to {\Cal N}$
with $H_{n+1}$ an extension of $H_n$. Since $n = 0$ is clear, we
assume that there is a homotopy $H_n \: \Lambda V_{(n)}^I \to
{\Cal N}$ with $H_n$ an  extension of $H_{n-1}$ and $f_n
\simeq_{H_n} 0$. Since $\Lambda V_{(n+1)} = \Lambda(V_{(n)}\oplus
V_{n+1})$ is an obstruction decomposition and $f_{n+1}\simeq 0$,
Proposition 2.8 shows that ${\Cal O}_{f_{n+1},0}^{H_n} =0$.  By
Proposition 2.5, $H_n$ can be extended to a homotopy $H_{n+1} \:
\Lambda V_{(n+1)}^I \to {\Cal N}$ with $f_{n+1} \simeq_{H_{n+1}} 0
\: \Lambda V_{(n+1)} \to {\Cal N}$.  This completes the inductive
step.  Finally, a homotopy $H \: \Lambda V^I \to {\Cal N}$ such
that $f \simeq_H 0 \: \Lambda V \to {\Cal N}$ is given by setting
$H|_{\Lambda V_{(n)}^I} = H_n$ for each $n$. \hfill$\square$
\enddemo

\noindent{\bf 3.2 Remarks.} (1) Theorem 3.1 gives a simple, direct
proof of the well-known fact that there are no non-trivial phantom
maps rationally \cite{McG, Th.\,3.20}.  We simply set $V_n =V^n $
so that $V_{(n)}=V^{(n)}$.  It would  be interesting to extend the
proof of Theorem 3.1 to the case of two maps $f,g \:\Lambda V \to
{\Cal N}$ such that $f_n \simeq g_n$ (cf. \cite{Su$_1$,
Lem.\,2.7}, \cite{D-R, Lem.\,2.4}).

\noindent{}(2) We observe that in Theorem 3.1 there is  freedom in
choosing the filtration on the minimal algebra $\Lambda V$.  In
particular, if a map $f \: \Lambda V \to {\Cal N}$ satisfies $f_n
\simeq 0$ for each $n$, for  {\it one}  filtration, then the same
must be true for {\it all} choices of filtration.
\smallskip

We next give a corollary which is essentially the contrapositive
of Theorem 3.1.  It is  this form which we will frequently use.

\proclaim{3.3 Corollary} Suppose $f\: {\Cal M}=\Lambda V \to {\Cal
N}$ represents a non-trivial homotopy class in $[{\Cal M},{\Cal
N}]$ and that $V = \oplus_{n \geq 0} V_n$ gives a filtration of
${\Cal M}$. Then there is  a map $f' \: {\Cal M} \to {\Cal N}$ and
an integer $n$ such that $f \simeq f' \: {\Cal M} \to {\Cal N}$,
$f'_{n-1} = 0 \: \Lambda V_{(n-1)} \to {\Cal N }$  and  the
obstruction ${\Cal O}_{f'_n,0} \: V_n \to H^*({\Cal N})$  is
non-zero.
\endproclaim

\demo{Proof} Suppose $f \not\simeq 0$.  By Theorem 3.1, there is
some $n$ for which $f_n \not\simeq 0 \: \Lambda V_{(n)} \to {\Cal
N}$. Let us choose the smallest such $n$. Then $f_{n-1} \simeq_H 0
\: \Lambda V_{(n-1)} \to {\Cal N}$ for some homotopy $H \: \Lambda
V_{(n-1)}^I \to {\Cal N}$.  Since the inclusion $\Lambda V_{(n-1)}
\to \Lambda V$ has the homotopy extension property, $H$ extends to
a homotopy from $\Lambda V^I$ to ${\Cal N}$ starting at $f$ and
ending at some map $f'$,  with $f'_{n-1} = f'|_{\Lambda V_{(n-1)}}
= 0$.   Thus the obstruction ${\Cal O}_{f'_n,0}$ is well-defined.
If ${\Cal O}_{f'_n,0} =0$, then  $f'_{n} \simeq 0 \: \Lambda
V_{(n)} \to {\Cal N}$.  Therefore $f_{n} \simeq 0 \: \Lambda
V_{(n)} \to {\Cal N}$ because $f_n \simeq f'_n$.  Since this is
not so by assumption, ${\Cal O}_{f'_n,0} \not= 0$. \hfill$\square$
\enddemo

\subhead{\bf \S4 A Conjecture of Copeland-Shar}\endsubhead
 Copeland and Shar have made the following
natural conjecture.

\proclaim{4.1 Conjecture} \cite{C-S, Conj.\,5.8} For rational
spaces $A$ and $B$, the set $[A,B]$ is either infinite or consists
of a single element.
\endproclaim

In this section, we establish Conjecture 4.1 under the
 hypothesis that either $A$ or $B$ is the rationalization of a universal space
(Theorem 4.6).  We give examples in Section 5 to show Conjecture
4.1 is false in general.  We begin with the definition of a
universal space.

\proclaim{\bf 4.2 Definition}\cite{M-O-T} A finite complex $X$ is
{\it universal\/} if for any ${\Bbb Q}$-equivalence $k \: U\to V$
and map $g \: X\to V$, there exists a ${\Bbb Q}$-equivalence $f \:
X\to X$ and a map $h \: X\to U$ such that $kh$ and $gf \: X \to V$
are homotopic.
\endproclaim

The class of universal spaces is large and includes many familiar
spaces.  For example, any formal space is universal by \cite{Sh}.
This means, in particular, that Eilenberg-MacLane spaces, Moore
spaces, $H$-spaces, co-$H$-spaces and K\"ahler manifolds are
universal.  In addition, all homogeneous spaces are universal,
including those that are not formal.
  Also, products and coproducts of
universal spaces are  universal.

Next we characterize universality in terms of minimal algebras.
\proclaim{\bf 4.3 Definition} If ${\Cal A}$ is a bigraded DG
algebra with each element of positive total degree, then ${\Cal
A}$ is called a {\it totally positive\/} bigraded DG algebra. A
{\it positive weight decomposition} on a graded algebra ${\Cal A}$
is a vector space decomposition ${\Cal A} = \oplus_{r \geq
1}\,{\Cal A}(r)$, such that ${\Cal A}(r)\cdot{\Cal A}(s) \subseteq
{\Cal A}(r+s)$ and $d\: {\Cal A}(r) \to {\Cal A}(r)$ (cf.
\cite{D-R}).  Elements of ${\Cal A}(r)$ are said to have weight
$r$.
\endproclaim

Given a positive weight decomposition of a minimal algebra ${\Cal
M}$, one defines a
 second grading on ${\Cal M}$ by  assigning to each homogeneous weight
generator $x$ a second degree  equal to $\text{weight}(x) - |x|$.
Then ${\Cal M}$ is a totally positive bigraded minimal algebra.
Conversely, given a totally positive bigraded minimal algebra
${\Cal M}$, one defines a positive weight decomposition on ${\Cal
M}$ by setting the weight equal to the total degree.

\proclaim{\bf 4.4 Definition} A minimal algebra ${\Cal M} =
\Lambda V$ is called {\it universal\/} if ${\Cal M}$ admits a
totally positive bigraded minimal algebra structure.
Equivalently, ${\Cal M}$ is universal if it admits a positive
weight decomposition.
\endproclaim

It is known that a space is  universal (Definition 4.2) if and
only if its Sullivan minimal model is universal (Definition 4.4)
\cite{Sc, Th.\,1}.

If ${\Cal M}=\Lambda V$ is a universal minimal algebra with $V =
\oplus_i V_i$ a totally positive second grading, then the second
grading on ${\Cal M}$ gives a decomposition $H^n({\Cal M}) \cong
\oplus _{i>-n} H_i^n({\Cal M})$.

\proclaim{4.5 Lemma} Let ${\Cal M}=\Lambda V$ be a universal
minimal algebra with $V = \oplus_i V_i$ a totally positive second
grading. Then for each $\lambda \not=0 \in {\Bbb Q}$, ${\Cal M}$
admits an automorphism $\phi_{\lambda} \: {\Cal M} \to {\Cal M}$
such that $\phi_{\lambda}(v) = \lambda^{|v|+i} v$ for $v \in V_i$.
Furthermore, for each element $x \in H_i^n({\Cal M})$, we have
$\phi_{\lambda}^*(x) = \lambda^{n+i} x$.
\endproclaim

\demo{Proof} Straightforward. \hfill$\square$
\enddemo

The existence of the self-maps and their induced `grading-like'
automorphisms on cohomology in Lemma 4.5 is a key property of
universal minimal algebras which we use in Theorem 4.6.

Since the homotopy category of minimal algebras is equivalent to
the homotopy category of rational spaces, the following result
establishes Conjecture 4.1 with an additional hypothesis.

\proclaim{4.6 Theorem} If ${\Cal M}$ and ${\Cal N}$ are minimal
algebras and ${\Cal M}$ or ${\Cal N}$ is universal, then the set
$[{\Cal M},{\Cal N}]$ is either infinite or consists of a single
element.
\endproclaim

\demo{Proof} We assume that $[{\Cal M},{\Cal N}]$ is not trivial
and show that $[{\Cal M},{\Cal N}]$ is infinite.  Choose the
degree filtration of ${\Cal M} = \Lambda V$.  By Corollary 3.3,
there is some $n$ and a map $f \: {\Cal M }\to{\Cal N}$, such that
$f_{n-1} = 0 \: \Lambda V^{(n-1)} \to {\Cal N}$ and ${\Cal
O}_{f_n,0} \: V^n \to H^n({\Cal N})$ is non-zero.

Now suppose that ${\Cal N}$ is universal. As in Lemma 4.5, we have
the self-maps $\phi_{\lambda} \: {\Cal N} \to {\Cal N}$.  We
claim, that for $\lambda \not= 0, \pm1$, the compositions
$\phi^i_{\lambda}f = \phi_{\lambda} \circ \cdots \circ
\phi_{\lambda } \circ f$ for each $i\ge 0$ are homotopically
distinct.  The restriction of $\phi ^i_{\lambda}f$ to $\Lambda
V^{(n)}$ gives $(\phi^i_{\lambda}f)_n = \phi^i_{\lambda}f_n$.
This is zero on $\Lambda V^{(n-1)}$, and so we have an obstruction
${\Cal O}_{\phi^i_{\lambda}f_n,\phi^j_{\lambda}f_n}$ for each $i$
and $j$.  We detect that this is non-zero as follows: Since ${\Cal
N}$ is universal, then $H^*({\Cal N})$ decomposes as $\oplus _r
H^*_r({\Cal N})$ and the projection onto some summand, say $p \:
H^n({\Cal N}) \to H^n_r({\Cal N})$, gives a non-trivial
composition $p{\Cal O}_{f_n,0} \: V^n \to H^n_r({\Cal N})$. Then
by Remark 2.7 and Lemma 4.5,
$$\align p{\Cal O}_{\phi^i_{\lambda}f_n,\phi^j_{\lambda}f_n} &=
p{\Cal O}_{\phi^i_{\lambda}f_n,0} - p{\Cal
O}_{\phi^j_{\lambda}f_n,0} \\ &= p(\phi^i_{\lambda})^*{\Cal
O}_{f_n,0} - p(\phi^j_{\lambda})^*{\Cal O}_{f_n,0} \\ &=
(\lambda^{i(n+r)} - \lambda^{j(n+r)})p{\Cal O}_{f_n,0}\ ,\\
\endalign
$$ which is non-zero for $i \not= j$.  It follows that ${\Cal
O}_{\phi^i_{\lambda}f_n,\phi^j_{\lambda}f_n} \not= 0$.  Hence from
Proposition 2.6, $\phi^i_{\lambda}f$ and $\phi^j_{\lambda}f$ are
not homotopic for $i \not= j$ and so $[{\Cal M }, {\Cal N}]$ is
infinite.

Now suppose that ${\Cal M}$ is universal. Then we have self-maps
$\phi_{\lambda} \: {\Cal M} \to {\Cal M}$ as in Lemma 4.5.  Here
we claim that for fixed $\lambda \not= 0, \pm1$, the compositions
$f\phi^i_{\lambda}$ for each $i$ represent distinct classes in
$[{\Cal M }, {\Cal N}]$. The restriction to $\Lambda V^{(n)}$ of
each such composition, $(f\phi^i_{\lambda})_n =
f_n\phi^i_{\lambda}$ is zero on $\Lambda V^{(n-1)}$ and so we have
an obstruction ${\Cal
O}_{f_n\phi^i_{\lambda},f_n\phi^j_{\lambda}}$ for each $i$ and
$j$.  We detect that this is non-zero as follows: Choose some
homogeneous second degree summand $V_k^n$ of $V^n$ with inclusion
$q \: V_k^n \to V^n$, such that the composition ${\Cal O}_{f_n,0}q
\: V_k^n \to H^n({\Cal N})$ is non-zero. Since $\phi_{\lambda}(v)
= \lambda^{n+k}v$ for $v \in V_k^n$,  we have that
$$ \align {\Cal O}_{f_n\phi^i_{\lambda},f_n\phi^j_{\lambda}}q &=
{\Cal O}_{f_n\phi^i_{\lambda},0}q - {\Cal
O}_{f_n\phi^j_{\lambda},0}q \\ &= (\lambda^{i(n+k)} -
\lambda^{j(n+k)}){\Cal O}_{f_n,0}q\ .\\
\endalign
$$
Since this is non-zero for $i \not= j$, it follows that ${\Cal
O}_{f_n\phi^i_{\lambda},f_n\phi^j_{\lambda}} \not= 0$. Hence from
Proposition 2.6, $f\phi^i_{\lambda}$ and $f\phi^j_{\lambda}$ are
not homotopic for $i \not= j$. Thus  $[{\Cal M }, {\Cal N}]$ is
infinite.
\enddemo

\subhead{\bf \S5 Examples}\endsubhead Theorem 4.6 establishes a
special case of Conjecture 4.1 which  holds for a large class of
spaces.  However, Conjecture 4.1 is not true without some
additional hypothesis, as the  examples below illustrate. Any
non-trivial minimal algebra admits the identity $\iota$ as a
self-homomorphism which is not homotopic to the zero homomorphism.
So if ${\Cal M}$ is any universal minimal algebra, then $[{\Cal
M}, {\Cal M}]$ is infinite by Theorem 4.6.  By way of contrast, we
next give some examples of non-universal minimal algebras that
have finitely many homotopy classes of
self-maps.\footnote"$^{({\dagger})}$"{These examples are
modifications of one of Halperin and Oprea (private
communication).}

\smallskip

\noindent{\bf 5.1 Example.}  Let ${\Cal M} = \Lambda(x_1, x_2,
y_1, y_2, y_3, z)$ with degrees of the generators given by $|x_1|
= 18$, $|x_2| = 22$, $|y_1| = 75$, $|y_2| = 79$, $|y_3| = 83$ and
$|z| = 197$. The differential is defined as follows:
$$\alignat3 dx_1 &= 0 &\qquad dy_1 &= x_1^3 x_2 &\qquad dz &= (y_1
x_2 - x_1 y_2)(y_2 x_2 - x_1 y_3) + x_1^{11} + x_2^{9}\\ dx_2 &= 0
&\qquad dy_2 &= x_1^2 x_2^2 &\qquad &= y_1 y_2 x_2^2 - y_1 y_3 x_1
x_2 + y_2 y_3 x_1^2  + x_1^{11} + x_2^{9}\\ & &\qquad dy_3 &= x_1
x_2^3 & &
\endalignat$$
It is easily verified that $d^2 = 0$ and so ${\Cal M}$ is a
minimal algebra.  Now we show that $[{\Cal M}, {\Cal M}]$ consists
of exactly 2 elements.

Let $f \: {\Cal M} \to {\Cal M}$ be any map.  Because of the
choice of degrees of the generators, there exist constants $a_1,
a_2, b_1, b_2, b_3,c \in {\Bbb Q}$ such that
$$\alignat3 f(x_1) &= \ a_1\,x_1 &\qquad f(y_1) &= \ b_1\,y_1
&\qquad f(z)   &= \ c\,z\\ f(x_2) &= \ a_2\,x_2 &\qquad f(y_2) &=
\ b_2\,y_2 & &\\ & &\qquad f(y_3) &= \ b_3\,y_3 & &
\endalignat$$
From $df(y_j) = fd(y_j)$, for $j = 1$, $2$ and $3$, we obtain the
following:
$$b_1 = a_1^3 a_2, \qquad b_2 = a_1^2 a_2^2 \quad \text{and} \quad
b_3 = a_1 a_2^3\ .\tag1$$
Then $df(z) = fd(z)$ implies that $a_1^{11} = c = a_2^{9}$ and $c
= b_1 b_2 a_2^2$.  From (*) it follows that $c = b_1b_2a_2^2  =
a_1^5a_2^5$.  We note that $a_1 = 0$ if and only if $a_2 = 0$.  In
this case, $b_1 = b_2 = b_3 = c = 0$ and so $f=0$.  Now suppose
that $a_1 \not= 0$, so that $a_2 \not= 0$. From $a_1^{11} = c =
a_1^5a_2^5$ we obtain
$$ a_1^6 = a_2^5 \tag2$$
and from $a_2^{9} = c = a_1^5a_2^5$ we obtain
$$ a_1^5 = a_2^4\ . \tag3$$
Dividing (2) by (3), we get $a_1 = a_2$.  From (2) we obtain $a_1
= 1$, and so $a_2 = b_1 = b_2 = b_3 = c = 1$ also.   In summary,
if $a_1 \not= 0$, then $f = \iota$.  This shows that the set of
homotopy classes of maps $[{\Cal M}, {\Cal M}]$ contains at most
two homotopy classes of maps, namely $[0]$ and $[\iota]$.  On the
other hand, these must be distinct classes, since ${\Cal M}$ is
not acyclic. Thus $[{\Cal M}, {\Cal M}]$ contains exactly two
elements.

\smallskip

\noindent{\bf 5.2 Example.} Let ${\Cal M} = \Lambda(x_1, x_2, y_1,
y_2, y_3, z)$ with degrees of these generators given by $|x_1| =
8$, $|x_2| = 10$, $|y_1| = 33$, $|y_2| = 35$, $|y_3| = 37$ and
$|z| = 119$. The differential is defined as follows:
$$\alignat3 dx_1 &= 0 &\qquad dy_1 &= x_1^3 x_2 &\qquad dz &=
x_1^4 (y_1 x_2 - x_1 y_2)(y_2 x_2 - x_1 y_3) + x_1^{15} +
x_2^{12}\\ dx_2 &= 0 &\qquad dy_2 &= x_1^2 x_2^2 &\qquad &=
y_1y_2x_1^4x_2^2 - y_1y_3x_1^5x_2 + y_2 y_3 x_1^6 + x_1^{15} +
x_2^{12}\\ & &\qquad dy_3 &= x_1 x_2^3 & &
\endalignat$$
Once again we show that $[{\Cal M}, {\Cal M}]$ consists of exactly
2 elements.

Let $f \: {\Cal M} \to {\Cal M}$ be any homomorphism.  For
constants $a_1, a_2, b_1, b_2, b_3 \in {\Bbb Q}$, we have
$$\alignat2 f(x_1) &= \ a_1\,x_1 &\qquad f(y_1) &= \ b_1\,y_1 \\
f(x_2) &= \ a_2\,x_2 &\qquad f(y_2) &= \ b_2\,y_2 \\ & &\qquad
f(y_3) &= \ b_3\,y_3
\endalignat$$
From $df(y_j) = fd(y_j)$, for $j = 1$, $2$ and $3$, we obtain the
following:
$$b_1 = a_1^3 a_2, \qquad b_2 = a_1^2 a_2^2 \quad \text{and} \quad
b_3 = a_1 a_2^3\ .\tag1$$
Also, since ${\Cal M}^{119} =\Lambda ( z,\, y_1 x_1^2 x_2^7,\, y_1
x_1^7 x_2^3,\, y_2 x_1^3 x_2^6,\,  y_2 x_1^8 x_2^2,\, y_3 x_1^4
x_2^5,\, y_3 x_1^9 x_2 )$, there exist $c, \lambda_1, \lambda_2,
\lambda_3, \nu_1, \nu_2,\nu_3$ $\in {\Bbb Q}$ such that
$$\align f(z)   = \ c\,z \ +\ &\lambda_1\,y_1\,x_1^2\,x_2^7 \;+\;
\lambda_2\,y_2\,x_1^3\,x_2^6\;+\;\lambda_3\,y_3\,x_1^4\,x_2^5\\
&\;+\;\nu_1\,y_1\,x_1^7\,x_2^3 \;+\;
\nu_2\,y_2\,x_1^8\,x_2^2\;+\;\nu_3\,y_3\,x_1^9\,x_2\ .\\
\endalign$$
Then $df(z) = fd(z)$ implies that $a_1^{15} = c = a_2^{12}$, $c =
b_1 b_2 a_1^4a_2^2$ and that $\lambda_1 + \lambda_2 + \lambda_3 =
0$ and $\nu_1 + \nu_2 + \nu_3 = 0$.  From (1) it follows that $c =
a_1^9a_2^5$. If $a_1 = 0$, then $a_2 = b_1 = b_2 = b_3 = c = 0$.
So suppose that $a_1 \not= 0$, so that $a_2 \not= 0$.  From
$a_1^{15} = c = a_1^9a_2^5$ we obtain
$$ a_1^6 = a_2^5 \tag2$$
and from $a_2^{12} = c = a_1^9a_2^5$ we obtain
$$ a_1^9 = a_2^7\ . \tag3$$
Dividing (3) by (2), we get
$$ a_1^3 = a_2^2 \ .\tag4$$
Squaring (4) and dividing (2) by the result gives $a_2 = 1$. From
(4) it follows that $a_1 = 1$ also.  In summary, there are two
distinct cases:

Case 1: $a_1 = a_2 = b_1 = b_2 = b_3 = c = 0$.

Case 2: $a_1 = a_2 = b_1 = b_2 = b_3 = c = 1$.

\noindent{}Next notice that in either case $\lambda_1 + \lambda_2
+ \lambda_3 = 0$ and $\nu_1 + \nu_2 + \nu_3 = 0$ imply that
$$\align f(z) &= \ c\,z \ -\;(\lambda_2 +
\lambda_3)\,y_1\,x_1^2\,x_2^7
\;+\;\lambda_2\,y_2\,x_1^3\,x_2^6\;+\;\lambda_3\,y_3\,x_1^4\,x_2^5\\
&\qquad\qquad \;-\;(\nu_2 + \nu_3)\,y_1\,x_1^7\,x_2^3
\;+\;\nu_2\,y_2\,x_1^8\,x_2^2\;+\;\nu_3\,y_3\,x_1^9\,x_2\\ &= cz +
\lambda_2 x_2^5 (- y_1 x_1^2 x_2^2 + x_1^3x_2y_2)
      + \lambda_3 x_1x_2^4 (- y_1 x_1 x_2^3 + x_1^3x_2y_3)\\
&\qquad\qquad + \nu_2 x_1^5x_2 (- y_1 x_1^2 x_2^2 + x_1^3x_2y_2)
      + \nu_3 x_1^6 (- y_1 x_1 x_2^3 + x_1^3x_2y_3)\\
&= cz + \lambda_2 x_2^5 d(y_1y_2) + \lambda_3 x_1x_2^4 d(y_1y_3)
      + \nu_2 x_1^5x_2 d(y_1y_2) + \nu_3 x_1^6 d(y_1y_3) \\
&= cz + d(\lambda_2x_2^5y_1y_2 + \lambda_3x_1x_2^4y_1y_3
      + \nu_2x_1^5x_2y_1y_2 + \nu_3x_1^6y_1y_3)\ .
\endalign
$$
Now write $V = V_0 \oplus V_1$ with $V_0 = \langle x_1, x_2, y_1,
y_2, y_3 \rangle$ and $V_1 = \langle z \rangle$. This gives an
obstruction decomposition of ${\Cal M} = \Lambda V$ and we can use
the results of Section 2. Consider any homomorphism $f \: {\Cal M}
\to {\Cal M}$ in Case 1 above.  Then $f|_{\Lambda V_0} = 0$. Let
$H \: \Lambda V_0^I \to {\Cal M}$ be the zero homotopy between $f$
and 0.  Then
$${\Cal O}^H_{f,0}(z) = [f(z)] = [d(\lambda_2x_2^5y_1y_2 +
\lambda_3x_1x_2^4y_1y_3
      + \nu_2x_1^5x_2y_1y_2 + \nu_3x_1^6y_1y_3)] = 0.$$
By Proposition 2.5, $f\simeq 0$. Likewise, for any homomorphism
$f$ in Case 2, we have that $f\simeq \iota$. This shows that the
set of homotopy classes of maps $[{\Cal M}, {\Cal M}]$ contains at
most two elements, namely $[0]$ and $[\iota]$.  On the other hand,
these must be distinct classes since ${\Cal M}$ is not acyclic.
Thus $[{\Cal M}, {\Cal M}]$ contains exactly two elements. We
contrast this example with the previous one in Remarks 5.4.

\smallskip

\noindent{\bf 5.3 Example.} Let ${\Cal N} = \Lambda(x_1, x_2, y_1,
y_2, y_3, z)$ with $|x_1| = 10$, $|x_2| = 12$, $|y_1| = 41$,
$|y_2| = 43$, $|y_3| = 45$ and $|z| = 119$. The differential is as
follows:
$$\alignat3 dx_1 &= 0 &\qquad dy_1 &= x_1^3 x_2 &\qquad dz &=
x_2(y_1 x_2 - x_1 y_2)(y_2 x_2 - x_1 y_3) + x_1^{12} + x_2^{10}\\
dx_2 &= 0 &\qquad dy_2 &= x_1^2 x_2^2 &\qquad &= y_1y_2x_2^3 -
y_1y_3 x_1x_2^2 + y_2y_3x_1^2x_2 + x_1^{12} + x_2^{10}\\ & &\qquad
dy_3 &= x_1 x_2^3 & &
\endalignat$$
Since ${\Cal N}^{119} = \langle z,\, y_1 x_1^3 x_2^4,\,  y_2 x_1^4
x_2^3,\, y_3 x_1^5 x_2^2 \rangle$, the possibilities for
self-homomorphisms are again not so restricted as in Example 5.1.
By an argument as in Example 5.2, one shows that $[{\Cal N}, {\Cal
N}]$ consists of exactly three elements.  We omit the details, but
assert that representatives of the three classes are the trivial
map 0, the identity $\iota$ and an involution $f$ given on
generators by $f(x_1) = x_1$, $f(x_2) = -x_2$, $f(y_1) = -y_1$,
$f(y_2) = y_2$, $f(y_3) = -y_3$ and $f(z) = z$.  It is easily seen
that none of these are homotopic by considering their induced
homomorphisms on cohomology.

\smallskip

\noindent{\bf 5.4 Remarks.}  These examples are interesting for a
number of reasons.  Example 5.1 is notable as the minimal algebra
itself admits only two self-homomorphisms,  the zero homomorphism
and the identity homomorphism. In general we  have to take into
account homotopy of homomorphisms as in Examples 5.2 and 5.3.

Recall that a minimal algebra $\Lambda V$ is called {\it
elliptic\/} if both $V$ and $H^*(\Lambda V)$ are
finite-dimensional \cite{Fe}. Example 5.1 has infinite-dimensional
cohomology and so is not elliptic. Examples 5.2 and 5.3 are
elliptic, however.  To see this for Example 5.2 we argue as
follows:   First note that
$$d(zx_2 - y_1y_2y_3x_1^3 - y_1x_1^{12}) = x_2^{13}.$$
Thus the class $[x_2] \in H^*({\Cal M})$ is nilpotent.  Now
consider the quotient minimal algebra $\overline{\Cal M} =
\Lambda(x_1, y_1, y_2, y_3, z)$, obtained by setting $x_2$ equal
to zero.  The differential $\overline d$ of $\overline{\Cal M}$ is
zero on all generators except $z$ and $\overline d(z) =
y_2y_3x_1^6 + x_1^{15}$.  We argue that the class $[x_1] \in
H^*(\overline{\Cal M})$ is nilpotent:  Re-order the generators of
$\overline{\Cal M}$ as $\Lambda(y_1, y_2, y_3, x_1, z)$.  An easy
application of Halperin's criterion \cite{Ha$_1$, Prop.\,1} now
shows $\overline{\Cal M}$ has finite-dimensional cohomology, and
hence $[x_1]$ is nilpotent in $H^*(\overline{\Cal M})$.  We
continue in this way by setting $x_1 =0$ and have $\Lambda(y_1,
y_2, y_3, z)$ as the next quotient.  Since all remaining
generators are of odd degree,  ${\Cal M}$ has finite-dimensional
cohomology and so is elliptic. A similar argument shows that
Example 5.3 is elliptic. Since elliptic spaces satisfy Poincar\'e
duality, it follows from \cite{Bar} that Examples 5.2 and 5.3 can
be taken as the minimal models of smooth manifolds.  Thus these
examples exhibit certain nice features despite the fact that their
generators are of high degrees.

Another interesting aspect of these examples concerns their group
of self-homo- topy equivalences.  For Examples 5.1 and 5.2 this
group is trivial, consisting of the class of $\iota$.  Because of
the equivalence between minimal algebras and rational spaces, this
provides examples of rational spaces with trivial group of
self-homotopy equivalences.  These complement the known examples
of such spaces (\cite{Ka$_1$} and \cite{Ka$_2$}), none of which is
rational.  Indeed, Example 5.2  gives an example of a
finite-dimensional (but not finite) CW-complex, which is a
rational space, with trivial group of self-homotopy equivalences.
For Example 5.3, the group of self-homotopy equivalences is ${\Bbb
Z}_2$. This is interesting, if for no other reason than to
illustrate the rather surprising appearance of finite groups in
rational homotopy theory! These examples raise questions about the
group of self-homotopy equivalences, such as the following:  Which
finite groups can be realized as the group of self-homotopy
equivalences of a rational space?

\enddocument